\documentclass[12pt]{article}
\usepackage{amssymb,amsmath}
\usepackage{cases}
\usepackage{amsfonts}
\usepackage{color,xcolor}
\usepackage[left=2.0cm,right=2.0cm,top=2.0cm,bottom=2.0cm]{geometry}
\usepackage[colorlinks,citecolor=blue,urlcolor=blue]{hyperref}
\usepackage[utf8]{inputenc}

\newtheorem{theorem}{Theorem}[section]

\newtheorem{lemma}{Lemma}[section]
\newtheorem{proposition}{Proposition}[section]

\newtheorem{definition}{Definition}[section]
\newtheorem{remark}{Remark}[section]

\newcommand{\bal}{\begin{align}}
\newcommand{\bbal}{\begin{align*}}
\newcommand{\beq}{\begin{equation}}
\newcommand{\eeq}{\end{equation}}
\newcommand{\bca}{\begin{cases}}
\newcommand{\eca}{\end{cases}}

\newcommand{\pa}{\partial}
\newcommand{\fr}{\frac}

\newcommand{\dd}{\mathrm{d}}

\newcommand{\R}{\mathbb{R}}

\newcommand{\les}{\lesssim}

\newcommand{\bi}{\Big}

\begin{document}
\title{Well-posedness and Continuity Properties of the Fornberg-Whitham equation in Besov space $B^1_{\infty,1}(\R)$}

\author{Guorong Qu$^{1}$,  Xing Wu$^{2,}$\footnote{E-mail: guorongqu@163.com; ny2008wx@163.com(Corresponding author); yuxiao5726@163.com} , Y. Xiao$^{2}$\\
\small $^1$ School of Tourism Data, Guilin Tourism University, Guilin 541006, China\\
\small $^2$ \small  College of Information and Management Science,
Henan Agricultural University,\\
\small Zhengzhou, Henan, 450046, China\\}

\date{}

\maketitle\noindent{\hrulefill}

{\bf Abstract:} For the Fornberg-Whitham equation, the local well-posedness in the critical Besov space $B_{p, 1}^{1+\frac{1}{p}}(\mathbb{R})$ with $1\leq p <\infty$ has been studied in \cite{Guo 2023}(Guo,  Nonlinear Anal. RWA., 2023). However, for the endpoint case $p=\infty$, whether it is locally well-posed or ill-posed in $B_{\infty, 1}^{1}(\mathbb{R})$ is still unknown. In this paper, we prove that the Fornberg-Whitham equation is well-posed in the critical Besov space $B_{\infty, 1}^{1}(\mathbb{R})$ with solutions depending
continuously on initial data, which is different from that of the Camassa-Holm equation \cite{Guo 2022}(Guo et al.,  J. Differ. Equ., 2022). In addition, we show that this dependence is
sharp by showing that the solution map is not uniformly continuous on the initial data.


\vskip0mm\noindent{\hrulefill}

\section{Introduction}\label{sec1}
In this paper, we consider the Cauchy problem of the  Fornberg-Whitham (FW) equation
\begin{align}\label{eq1}
\begin{cases}
u_{xxt}-u_t+\frac{9}{2}u_xu_{xx}+\frac{3}{2}uu_{xxx}-\frac{3}{2}uu_x+u_x=0, &\quad (t,x)\in \mathbb{R}^+\times\mathbb{R},\\
u(0,x)=u_0(x), &\quad x\in \mathbb{R},
\end{cases}
\end{align}
which was first introduced by Whitham  and Fornberg \cite{Whitham 1967, Whitham 1978}  to study the nonlinear wave breaking phenomena. Using the Green function $G(x)=\frac{1}{2}e^{-|x|}$ and the identity $(1-\partial_{xx})^{-1}f=G*f$ for all $f\in L^2(\mathbb{R})$, we
can rewrite (\ref{eq1}) in the following non-local form
\begin{align}\label{eq2}
\begin{cases}
u_t+\frac32uu_x=\partial_x(1-\partial_{xx})^{-1}u=\partial_xG*u, &\quad (t,x)\in \R^+\times\R,\\
u(0,x)=u_0(x), &\quad x\in \R.
\end{cases}
\end{align}
In contrast to the famous Korteweg-de Vries (KdV) equation \cite{Korteweg 1985}
 \begin{align*}
   u_t+6uu_x=-\partial_{xxx}u,
 \end{align*}
 which can not describe the wave breaking phenomena, and  the classical Camassa-Holm (CH) equation \cite{Camassa 1993, Constantin 1998, Constantin 2009}
  \begin{align*}
   u_t+uu_x=-\partial_x(1-\partial_{xx})^{-1}(u^2+\frac{1}{2}u_x^2),
 \end{align*}
the FW equation not only admits solitary traveling wave solutions like
the KdV equation, but also has peaked traveling wave solutions
as the CH equation which are of the form $u(t, x)=\frac{8}{9}e^{-\frac{1}{2}|x-\frac{4}{3}t|}\cite{Whitham 1978}.$ Moreover, unlike the KdV and CH equations that are integrable and have infinitely many conserved quantities, the FW equation is not integrable. Although the FW equation is presented in a simple form,
the only useful conservation law we know so far is $\|u\|_{L^2}.$ Therefore, the analysis of the FW
equation would be somewhat more difficult.

The local posedness (existence, uniqueness and continuous Dependence) was first established by Holmes\cite{Holmes 2016} in the Sobolev spaces $H^s$ with $s>\frac{3}{2}$ on the circle, and  later extended by Holmes-Thompson \cite{Holmes 2017} to the Besov space $B_{2, r}^s$ with $s>\frac{3}{2}, 1<r<\infty$ or $s=\frac{3}{2}, r=1.$  Recently, Guo\cite{Guo 2023} investigated the  FW equation in  a general Besov space, and obtained the local well-posedness in $B_{p, r}^s(\mathbb{R})$, $s>1+\frac{1}{p}$, $1\leq p\leq \infty,$ $1\leq r< \infty$ and critical Besov spaces $B_{p, 1}^{1+\frac{1}{p}}(\mathbb{R})$ with $1\leq p <\infty.$ Li et al.\cite{Li 2023} proved that  the  FW equation  is ill-posed in Besov space $B_{p, 1}^s(\mathbb{R})$ with $1\leq s<1+\frac{1}{p}$ and $2\leq p<\infty$ due to the norm inflation. For the endpoint case $p=\infty$, whether it is locally well-posed or ill-posed in $B_{\infty, 1}^{1}(\mathbb{R})$ is unknown. Guo et al. \cite{Guo 2022} proved the ill-posedness for the CH equation in $B_{\infty, 1}^{1}(\mathbb{R})$ due to the norm inflation. However, motivated by \cite{Guo 2019}, by fully utilizing the structure of the equation and the internal properties of $u_0$ lying in  $B_{\infty, 1}^{1}(\mathbb{R})$ , we show the well-posedness for the FW equation in critical Besov spaces $B_{\infty, 1}^{1}(\mathbb{R})$.

 Our main result is stated as follows.
\begin{theorem}\label{the1.1} Let $u_0\in B_{\infty, 1}^{1}(\mathbb{R}).$ Then there exists a time $T>0$ such that  the Cauchy problem (\ref{eq2}) has a unique solution in $\mathcal{C}([0, T];B_{\infty, 1}^{1}(\mathbb{R}))\cap \mathcal{C}^1([0, T];B_{\infty, 0}^{1}(\mathbb{R}))$ and
 the data-to-solution map $u_0\rightarrow \mathcal{S}_t(u_0)$ is continuous from any bounded subset of  $B_{\infty, 1}^{1}(\mathbb{R})$ to $\mathcal{C}([0, T];B_{\infty, 1}^{1}(\mathbb{R}))$. That is, the FW equation is locally well-posed in $B_{\infty, 1}^{1}(\mathbb{R})$ in the sense of Hadamard.
\end{theorem}
\begin{remark}\label{rem1}
Since the CH equation is locally ill-posed in $B_{\infty, 1}^{1}(\mathbb{R})$, but for FW equation, we obtain the local well-posedness in $B_{\infty, 1}^{1}(\mathbb{R})$.   This interesting fact illustrates that there is a nature difference between the two equations.
\end{remark}

From our well-posedness result, we can further show that the dependence of the solutions on initial data in $B_{\infty, 1}^{1}(\mathbb{R})$ can not be better than continuous. More precisely, we shall prove the following result.
\begin{theorem}\label{the1.2} The data-to-solution map $u_0\rightarrow \mathcal{S}_t(u_0)$ for the FW equation, defined by the Cauchy
problem (\ref{eq2}),  is not uniformly continuous from any
bounded subset in  $B_{\infty, 1}^{1}(\mathbb{R})$ into $\mathcal{C}([0, T];B_{\infty, 1}^{1}(\mathbb{R}))$. More precisely,
there exist two sequences of solutions $f_n$ and $g_n$ such that
\begin{align*}
 \|f_n\|_{B_{\infty, 1}^{1}}\lesssim 1 \qquad \qquad  and  \qquad \qquad \lim_{n\rightarrow\infty}\|g_n\|_{B_{\infty, 1}^{1}}=0,
\end{align*}
but
\begin{eqnarray*}
  \liminf_{n\rightarrow \infty} \|\mathcal{S}_t(f_n+g_n)-\mathcal{S}_t(f_n)\|_{B_{\infty, 1}^1}\gtrsim t, \quad t\in [0, T_0],
        \end{eqnarray*}
with small positive time $T_0$ for $T_0\leq T$.
\end{theorem}

{\bf Notations}:  Given a Banach space $X$, we denote the norm of a function on $X$ by $\|\|_{X}$, and \begin{eqnarray*}
\|\cdot\|_{L_T^\infty(X)}=\sup_{0\leq t\leq T}\|\cdot\|_{X}.
\end{eqnarray*}
For $\mathbf{f}=(f_1, f_2,...,f_n)\in X$,
\begin{eqnarray*}
\|\mathbf{f}\|_{X}^2=\|f_1\|_{X}^2+\|f_2\|_{X}^2+...+\|f_n\|_{X}^2.
\end{eqnarray*}
The symbol
$A\lesssim B$ means that there is a uniform positive constant $C$ independent of $A$ and $B$ such that $A\leq CB$.

\section{Preliminaries}\label{sec2}
In this section, we  review the definition of Littlewood-Paley decomposition and nonhomogeneous Besov space, and then list some useful properties. For more details, the readers can refer to \cite{Bahouri 2011}.

\begin{definition} Let $\mathcal{B}:=\{\xi\in\mathbb{R}:|\xi|\leq 4/3\}$ and $\mathcal{C}:=\{\xi\in\mathbb{R}:3/4\leq|\xi|\leq 8/3\}.$
There exist two radial functions $\chi\in C_c^{\infty}(\mathcal{B})$ and $\varphi\in C_c^{\infty}(\mathcal{C})$ both taking values in $[0,1]$ such that
$$
\chi(\xi)+\sum_{j\geq0}\varphi(2^{-j}\xi)=1 \quad \forall \;  \xi\in \R.$$
For every $u\in \mathcal{S'}(\mathbb{R})$, the Littlewood-Paley dyadic blocks ${\Delta}_j$ are defined as follows
\begin{numcases}{\Delta_ju=}
0, & if $j\leq-2$;\nonumber\\
\chi(D)u=\mathcal{F}^{-1}(\chi \mathcal{F}u), & if $j=-1$;\nonumber\\
\varphi(2^{-j}D)u=\mathcal{F}^{-1}\big(\varphi(2^{-j}\cdot)\mathcal{F}u\big), & if $j\geq0$.\nonumber
\end{numcases}
The inhomogeneous low-frequency cut-off operator $S_j$ is defined by
$$
S_ju=\sum_{q=-1}^{j-1}{\Delta}_qu.
$$

\end{definition}
\begin{definition}
Let $s\in\mathbb{R}$ and $(p,r)\in[1, \infty]^2$. The nonhomogeneous Besov space $B^{s}_{p,r}(\R)$ is defined by
\begin{align*}
B^{s}_{p,r}(\R):=\Big\{f\in \mathcal{S}'(\R):\;\|f\|_{B^{s}_{p,r}(\mathbb{R})}<\infty\Big\},
\end{align*}
where
\begin{numcases}{\|f\|_{B^{s}_{p,r}(\mathbb{R})}=}
\left(\sum_{j\geq-1}2^{jrs}\|\Delta_jf\|^r_{L^p(\mathbb{R})}\right)^{\fr1r}, &if $1\leq r<\infty$,\nonumber\\
\sup_{j\geq-1}2^{js}\|\Delta_jf\|_{L^p(\mathbb{R})}, &if $r=\infty$.\nonumber
\end{numcases}
\end{definition}
\begin{lemma}(\cite{Bahouri 2011})\label{lem2.1} \;\;Let $s>0$ and $1\leq p\leq \infty$.\\
(1) $B^s_{p,1}(\R)\cap L^\infty(\R)$ is a Banach  algebra. Moreover, $B^{0}_{\infty,1}(\R)\hookrightarrow L^\infty(\R)\hookrightarrow B^{0}_{\infty,\infty}(\R)$.\\
(2) \begin{align*}
   \lim_{j\rightarrow \infty}\|S_ju-u\|_{B_{p, 1}^s}=0
    \end{align*}
(3)For any $u,v \in B^s_{p,1}(\R)\cap L^\infty(\R)$, we have
 \begin{align*}
&\|uv\|_{B^{s}_{p,1}}\leq C(\|u\|_{B^{s}_{p,1}}\|v\|_{L^\infty}+\|v\|_{B^{s}_{p,1}}\|u\|_{L^\infty}).
\end{align*}
We also have the following interpolation inequality
\begin{align*}
&\|u\|_{B^{1}_{\infty,1}}\leq C\|u\|_{B^{0}_{\infty,\infty}}^{\fr12}\|u\|_{B^{2}_{\infty,\infty}}^{\fr12}.
\end{align*}
(4) Let $m\in \mathbb{R}$ and $f$ be an $S^m-$ multiplier (i.e., $f: \mathbb{R}\rightarrow \mathbb{R}$ is smooth and satisfies that $\forall \alpha\in \mathbb{N}$, there exists a constant $\mathcal{C}_\alpha$ such that $|\partial^\alpha f(\xi)|\leq \mathcal{C}_\alpha(1+|\xi|)^{m-\alpha}$ for all $\xi \in \mathbb{R}$). Then the operator $f(D)$ is continuous from $B_{p, 1}^s(\mathbb{R})$ to $B_{p, 1}^{s-m}(\mathbb{R})$.
\end{lemma}

\begin{lemma}\label{lem2.2}(\cite{Bahouri 2011, Li 2017})
Let $1\leq p\leq \infty$, $\sigma\geq 1+\frac{1}{p}.$
There exists a constant $C=C(p,\sigma)$ such that for any solution to the
following linear transport equation:
\begin{equation*}
\partial_t f+v\partial_x f=g,\qquad
f|_{t=0} =f_0,
\end{equation*}
the following statements hold:
\begin{align*}
  \|f(t)\|_{B^{\sigma}_{p,1}}\leq \|f_0\|_{B^{\sigma}_{p,1}}+\int_0^t\|g(\tau)\|_{B^{\sigma}_{p,1}}d\tau+\int_0^tCV^{'}(\tau)\|f(\tau)\|_{B^{\sigma}_{p,1}}d\tau
\end{align*}
or
\begin{align*}
\sup_{s\in [0,t]}\|f(s)\|_{B^{\sigma}_{p,1}}\leq e^{CV(t)}\Big(\|f_0\|_{B^\sigma_{p,1}}
+\int^t_0e^{-CV(\tau)}\|g(\tau)\|_{B^{\sigma}_{p,1}}\dd \tau\Big),
\end{align*}

with
\begin{align*}
V(t)=
\int_0^t \|\partial_x v(s)\|_{B^{\sigma-1}_{p,1}}\dd s
\end{align*}
\end{lemma}
\section{Proof of Theorem \ref{the1.1}}
\setcounter{equation}{0}

In this section, we divide Proof of Theorem \ref{the1.1} into three steps.
\subsection{Existence and Uniqueness}
To prove the local existence, it is sufficient to establish  a priori estimate of the  solution in $B_{\infty, 1}^1$.

Using Lemma \ref{lem2.2}-Lemma \ref{lem2.1}, we have
\begin{align}\label{eq3.1}
\|u(t)\|_{B^{1}_{\infty,1}}&\leq e^{CV(t)}\Big(\|u_0\|_{B^{1}_{\infty,1}}+\int_0^te^{-CV(\tau)}\|\partial_x(1-\partial_{xx})^{-1}u\|_{B^{1}_{\infty,1}}\dd \tau\Big)\nonumber\\
&\leq e^{CV(t)}\Big(\|u_0\|_{B^{1}_{\infty,1}}+\int_0^te^{-V(\tau)}\|u(\tau)\|_{B^{1}_{\infty,1}}\dd \tau\Big),
\end{align}
where $V(t)=\int_0^t \|u(\tau)\|_{B^{1}_{\infty,1}}\dd\tau$. Then we obtain from \eqref{eq3.1} that
\bbal
F(t):=e^{-CV(t)}\|u(t)\|_{B^{1}_{\infty,1}}&\leq \|u_0\|_{B^{1}_{\infty,1}}+\int_0^tF(\tau)\dd \tau,
\end{align*}
which alongs with Gronwall's inequality yields
\begin{align}\label{eq3.2}
\|u(t)\|_{B^{1}_{\infty,1}}&\leq C\|u_0\|_{B^{1}_{\infty,1}}\exp\Big(C\int_0^t \|u(\tau)\|_{B^{1}_{\infty,1}}\dd\tau\Big).
\end{align}
Let
\begin{align*}
 A(t)=C\|u_0\|_{B^{1}_{\infty,1}}\exp\Big(C\int_0^t \|u(\tau)\|_{B^{1}_{\infty,1}}\dd\tau\Big)\;\; \mbox{with}\; \;A(0)=C\|u_0\|_{B^{1}_{\infty,1}},
\end{align*}
then from (\ref{eq3.2}), we have
\begin{align*}
\frac{\dd A(t)}{\dd t}\leq CA^2(t).
\end{align*}

Solving the above differential inequalities and combining with  \ref{eq3.2}, one has
\begin{align}\label{eq3.3}
\|u(t)\|_{B^{1}_{\infty,1}}&\leq \frac{C\|u_0\|_{B^{1}_{\infty,1}}}{1-Ct\|u_0\|_{B^{1}_{\infty,1}}}
\end{align}
Fix a $T > 0$ such that $T<\frac{1}{2C\|u_0\|_{B^{1}_{\infty,1}}}$, then we have
\begin{align}\label{eq3.4}
 \|u\|_{L^\infty_TB^{1}_{\infty,1}}&\leq C\|u_0\|_{B^{1}_{\infty,1}}, \qquad \forall t\in [0, T].
\end{align}
Therefore, the solution $u$ is  uniformly bounded in $L^\infty([0, T]; B^{1}_{\infty,1}).$

The existence follows the standard procedure, we omit the details. The uniqueness is the direct result of the following lemma. In fact, suppose that $u_1,u_2\in \mathcal{C}([0,T],B^{1}_{\infty,1})$ are two solutions of (\ref{eq2}) with the same initial data $u_0$, then
we have
$$\|u_1-u_2\|_{L^\infty}\leq C\|u_1(0)-u_2(0)\|_{L^\infty}=0,$$
which implies the uniqueness.
\begin{lemma}\label{ley1} Let $u,v\in \mathcal{C}([0,T],B^{1}_{\infty,1})$ be two solutions of (\ref{eq2}) associated with $u_0$ and $v_0$, respectively. Then we have the estimate for the difference $w=u-v$
\begin{align}
&\|w\|_{L^\infty}\leq C\|w_0\|_{L^\infty},\label{eq3.5}
\\&\|w\|_{B^1_{\infty,1}}\leq C\Big(\|w_0\|_{B^1_{\infty,1}}+\int^t_0\|v_x\|_{B^1_{\infty,1}}\|w\|_{L^\infty}\dd \tau\Big),\label{eq3.6}
\end{align}
where the constants $C$ depends on $T$ and initial norm $\|u_0,v_0\|_{B^1_{\infty,1}}$.
\end{lemma}
{\bf Proof.}\quad
It is easy to check that $w$ satisfys
\begin{align}\label{eq3.7}
\begin{cases}
\pa_tw+\frac32uw_x=-\frac32wv_x+\pa_x(1-\pa_{xx})^{-1}w, \\
w(0,x)=u_0(x)-v_0(x).
\end{cases}
\end{align}
Let $p\geq 2$, taking the inner product of (\ref{eq3.7}) with $p|w|^{p-2}w$, then one has
\begin{align*}
  \frac{\dd\|w\|^p_{L^p}}{\dd t}&=\frac{3}{2}\int u_x|w|^p\dd t-\frac{3}{2}\int p|w|^pv_x\dd t+\int_t(\partial_xG*w)p|w|^{p-2}w\dd t.
\end{align*}
Using the Young inequality and H\"{o}lder inequality, we obtain
\begin{align*}
  \frac{\dd\|w\|^p_{L^p}}{\dd t}&\leq C(\|u_x\|_{L^\infty}+p\|v_x\|_{L^\infty}+p)\|w\|^p_{L^p},
\end{align*}
thus
\begin{align*}
  \frac{\dd\|w\|_{L^p}}{\dd t}&\leq C(\|u_x\|_{L^\infty}+\|v_x\|_{L^\infty}+1)\|w\|_{L^p}.
\end{align*}
Using the embedding $B_{\infty, 1}^1(\mathbb{R})\hookrightarrow \mathcal{C}^{0, 1}(\mathbb{R})$, and integrating the above differential inequality with respect to time $t\in [0, T]$ yields
\begin{align*}
\|w\|_{L^p}
&\leq \|w_0\|_{L^p}\exp\Big(C\int_0^T\big(\|v\|_{B^1_{\infty,1}}+\|v\|_{B^1_{\infty,1}}+1\big)\dd s\Big)\leq C\|w_0\|_{L^p}.
\end{align*}
Let $p\rightarrow\infty$, we obtain the desired estimate (\ref{eq3.5}).

Applying Lemma Lemma \ref{lem2.2}-Lemma \ref{lem2.1} to  (\ref{eq2}) yields
\begin{align*}
\|w(t)\|_{B^1_{\infty,1}}-\|w_0\|_{B^1_{\infty,1}}&\les\int_0^t\Big( \|u\|_{B^1_{\infty,1}}\|w\|_{B^1_{\infty,1}}+\|wv_x,\pa_x(1-\pa_{xx})^{-1}w\|_{B^1_{\infty,1}}\Big)\mathrm{d}\tau\\
&\les \int_0^t(\|u\|_{B^1_{\infty,1}}+1)\|w\|_{B^1_{\infty,1}}+\|wv_x\|_{B^1_{\infty,1}}\mathrm{d}\tau\\
&\les\int_0^t\|w\|_{B^1_{\infty,1}}\|v_x\|_{L^\infty}+\|v_x\|_{B^1_{\infty,1}}\|w\|_{L^\infty}+(\|u\|_{B^1_{\infty,1}}+1)\|w\|_{B^1_{\infty,1}}\mathrm{d}\tau\\
&\les\int_0^t\|w\|_{B^1_{\infty,1}}\big(1+\|u, v\|_{B^1_{\infty,1}}\big)\mathrm{d}\tau+\int_0^t\|v_x\|_{B^1_{\infty,1}}\|w\|_{L^\infty}\mathrm{d}\tau.
\end{align*}
Gronwall's inequality yields  the desired (\ref{eq3.6}).

\subsection{Continuous Dependence}
Now we will prove that the solution of the FW equation in $B_{\infty, 1}^1(\mathbb{R})$ is continuously dependent on the initial data by \cite{Guo 2019}. The main difficulty lies in that the FW equation is of hyperbolic type. Precisely speaking, if $u,v\in \mathcal{C}([0,T],B^{1}_{\infty,1})$ are two solutions of \eqref{eq1} associated with $u_0$ and $v_0$, in view of \eqref{eq3.6}, we have to tackle with the term $\|v_x\|_{B^1_{\infty,1}}$. To bypass this, we can take $v=\mathcal{S}_{t}(S_Nu_0)$ as the solution to \eqref{eq1} with initial data $S_Nu_0$.

Letting $u=\mathcal{S}_{t}(u_0)$ and $v=\mathcal{S}_{t}(S_Nu_0)$, using Lemma \ref{lem2.2} and \eqref{eq3.5}, we have
\bbal
&\mathcal{S}_{t}(S_Nu_0)\|_{B^2_{\infty,1}}\leq C\|S_Nu_0\|_{B^2_{\infty,1}}\leq C2^N\|u_0\|_{B^1_{\infty,1}}
\end{align*}
and
\bbal
\|\mathcal{S}_{t}(S_Nu_0)-\mathcal{S}_{t}(u_0)\|_{L^\infty}&\leq C\|S_Nu_0-u_0\|_{L^\infty}
\\&\leq C\|S_Nu_0-u_0\|_{B^0_{\infty,1}}\\
&\leq C2^{-N}\|S_{N}u_0-u_0\|_{B^1_{\infty,1}},
\end{align*}
which combining with (\ref{eq3.6}) imply
\bbal
\|\mathcal{S}_{t}(S_Nu_0)-\mathcal{S}_{t}(u_0)\|_{B^1_{\infty,1}}&\leq C\|S_Nu_0-u_0\|_{B^1_{\infty,1}}.
\end{align*}
Then, for $u_0,\widetilde{u}_0\in B^1_{\infty,1}$, we have
\bbal
\|\mathcal{S}_{t}(u_0)-\mathcal{S}_{t}(\widetilde{u}_0)\|_{B^1_{\infty,1}}
&\leq \|\mathcal{S}_{t}(S_Nu_0)-\mathcal{S}_{t}(u_0)\|_{B^1_{\infty,1}}+\|\mathcal{S}_{t}(S_N\widetilde{u}_0)-\mathcal{S}_{t}(\widetilde{u}_0)\|_{B^1_{\infty,1}}\\
&~~+\|\mathcal{S}_{t}(S_Nu_0)-\mathcal{S}_{t}(S_N\widetilde{u}_0)\|_{B^1_{\infty,1}}
\\&\leq C\|S_Nu_0-u_0\|_{B^1_{\infty,1}}+C\|S_N\widetilde{u}_0-\widetilde{u}_0\|_{B^1_{\infty,1}}\\
&~~+C\|\mathcal{S}_{t}(S_Nu_0)-\mathcal{S}_{t}(S_N\widetilde{u}_0)\|_{B^1_{\infty,1}}
\\
&:=\mathcal{I}_1+\mathcal{I}_2+\mathcal{I}_3.
\end{align*}
Using the interpolation inequality in Lemma \ref{lem2.1}, we obtain
\bbal
\mathcal{I}_3&\leq C\|\mathcal{S}_{t}(S_Nu_0)-\mathcal{S}_{t}(S_N\widetilde{u}_0)\|^{\fr12}_{B^0_{\infty,\infty}}
\|\mathcal{S}_{t}(S_Nu_0)-\mathcal{S}_{t}(S_N\widetilde{u}_0)\|^{\fr12}_{B^2_{\infty,\infty}}\\
&\leq C\|\mathcal{S}_{t}(S_Nu_0)-\mathcal{S}_{t}(S_N\widetilde{u}_0)\|^{\fr12}_{L^\infty}
\|\mathcal{S}_{t}(S_Nu_0)-\mathcal{S}_{t}(S_N\widetilde{u}_0)\|^{\fr12}_{B^2_{\infty,1}}\\
&\leq C\|S_Nu_0-S_N\widetilde{u}_0\|^{\fr12}_{L^\infty}
\|\mathcal{S}_{t}(S_Nu_0)-\mathcal{S}_{t}(S_N\widetilde{u}_0)\|^{\fr12}_{B^2_{\infty,1}}\\
&\leq C2^{\frac N2}\|u_0-\widetilde{u}_0\|^{\frac12}_{B^1_{\infty,1}},
\end{align*}
which clearly implies
\bbal
\|\mathcal{S}_{t}(u_0)-\mathcal{S}_{t}
(\widetilde{u}_0)\|_{B^1_{\infty,1}}
&\lesssim \|S_Nu_0-u_0\|_{B^1_{\infty,1}}+
\|S_N\widetilde{u}_0-\widetilde{u}_0\|_{B^1_{\infty,1}}+2^{\frac N2}\|u_0-\widetilde{u}_0\|^{\frac12}_{B^1_{\infty,1}}
\\&\lesssim \|S_Nu_0-u_0\|_{B^1_{\infty,1}}+
\|S_N\widetilde{u}_0-\widetilde{u}_0\|_{B^1_{\infty,1}}+2^{\frac N2}\|u_0-\widetilde{u}_0\|^{\frac12}_{B^1_{\infty,1}}.
\end{align*}
Taking $N$ large enough, we can conclude that the   the data-to-solution map $u_0\rightarrow \mathcal{S}_t(u_0)$ is  continuously dependent on initial data.

\section{Proof of Theorem \ref{the1.2}}
\setcounter{equation}{0}

 Let $\hat{\phi}\in \mathcal{C}^\infty_0(\mathbb{R})$ be an even, real-valued and non-negative function on $\R$ and satisfy
\begin{numcases}{\hat{\phi}(x)=}
1, &if $|x|\leq \frac{1}{4}$,\nonumber\\
0, &if $|x|\geq \frac{1}{2}$.\nonumber
\end{numcases}
\begin{lemma}\label{lem4.1}
We define the high frequency function $f_n$ and the low frequency functions $g_n$ as follows
\bbal
&f_n=2^{-n}\phi(x)\sin \bi(\frac{17}{12}2^nx\bi),\\
&g_n=\frac{12}{17}2^{-n}\phi(x),\quad n\gg1.
\end{align*}
Then for any $\sigma\in\R$, we have
\bal
&\|f_n\|_{L^\infty}\leq C2^{-n}\phi(0)\quad\text{and}\quad\|g_n\|_{L^\infty}\leq C2^{-n}\phi(0),\label{y0}\\
&\|f_n\|_{B^{\sigma}_{\infty,1}}\leq C2^{(\sigma-1)n}\phi(0)\quad\text{and}\quad\|g_n\|_{B^\sigma_{\infty,1}}\leq C2^{-(n+\sigma)}\phi(0),\label{yz}\\
&\liminf_{n\rightarrow \infty}\|g_n\pa_xf_n\|_{B^{1}_{\infty,\infty}}\geq M_1,\label{yz3}
\end{align}
for some positive constants $C, M_1$.
\end{lemma}
{\bf Proof.}\quad We refer to see Lemma 3.2-Lemma 3.4 in \cite{Li 2020} for the proof with minor modifications.
\begin{proposition}\label{pro1}
Assume that $\|u_0\|_{B^{1}_{\infty,1}}\lesssim 1$. Under the assumptions of Theorem \ref{the1.1}, we have
\bal\label{et0}
\|\mathcal{S}_{t}(u_0)-u_0-t\mathbf{v}_0(u_0)\|_{B^{1}_{\infty,1}}\leq Ct^{2}\mathcal{E}(u_0),
\end{align}
here $\mathbf{v}_0(u_0):=\partial_x(1-\partial_{xx})^{-1}u_0-\frac{3}{2}u_0\pa_x u_0$ and
\bbal
\mathcal{E}(u_0)&:=1+\|u_0\|_{L^\infty}\big(\|u_0\|_{B^{2}_{\infty,1}}+
\|u_0\|_{L^\infty}
\|u_0\|_{B^{3}_{\infty,1}}\big).
\end{align*}
\end{proposition}
{\bf Proof.}\quad For simplicity, we denote $u(t)=\mathbf{S}_t(u_0)$.
By the Mean Value Theorem and (\ref{eq3.4}), we obtain
\bal\label{et1}
\|u(t)-u_0\|_{L^\infty}&\leq\int^t_0\|\pa_\tau u\|_{L^\infty} \dd\tau
\nonumber\\&\leq \int^t_0\frac{3}{2}\|u\pa_xu\|_{L^\infty}\dd \tau+\int^t_0\|\partial_xG*u\|_{L^\infty}\dd \tau
\nonumber\\&\leq C\int^t_0\|u(\tau)\|_{L^\infty}\big(\|u(\tau)\|_{C^{0,1}}+1\big)\dd \tau\nonumber\\
&\leq Ct\|u_0\|_{L^\infty}.
\end{align}

Using Lemma \ref{lem2.1}, \eqref{eq3.4} yield
\bal\label{et2}
\|u(t)-u_0\|_{B^{1}_{\infty,1}}
&\leq \int^t_0\|\pa_\tau u\|_{B^{1}_{\infty,1}} \dd\tau
\nonumber\\&\leq \int^t_0\|\partial_x(1-\partial_{xx})^{-1}u\|_{B^{1}_{\infty,1}} \dd\tau+ \int^t_0\frac{3}{2}\|u \pa_xu\|_{B^{1}_{\infty,1}} \dd\tau
\nonumber\\&\leq Ct\big(\|u\|_{L^\infty _T(B^{1}_{\infty,1})}+\|u\|^2_{L^\infty _T(B^{1}_{\infty,1})}+\|u\|_{L^\infty _T(L^\infty)}\|u_x\|_{L_T^\infty(B^{1}_{\infty, 1})}\big)
\nonumber\\&\leq Ct\big(1+\|u_0\|_{L^\infty}
\|u_0\|_{B^{2}_{\infty,1}}\big).
\end{align}
Similarly, we have
\bal\label{et3}
\|u(t)-u_0\|_{B^{2}_{\infty,1}}
&\leq \int^t_0\|\pa_\tau u\|_{B^{2}_{\infty,1}} \dd\tau
\nonumber\\&\leq \int^t_0\|\partial_x(1-\partial_{xx})^{-1}u\|_{B^{2}_{\infty,1}} \dd\tau+ \int^t_0\frac{3}{2}\|u \pa_xu\|_{B^{2}_{\infty,1}} \dd\tau
\nonumber\\&\leq Ct\big(\|u\|_{L^\infty_T(B^{2}_{\infty,1})}
+\|u\|_{L^\infty_T(L^\infty)}\|u\|_{L^\infty_T(B^{3}_{\infty,1})}\big)
\nonumber\\&\leq Ct\big(\|u_0\|_{B_{\infty, 1}^2}
+\|u_0\|_{L^\infty}\|u_0\|_{B^{3}_{\infty,1}}\big).
\end{align}
Using the Mean Value Theorem and Lemma \ref{lem2.1} once again, we obtain that
\bal\label{et4}
\|u(t)-u_0-t\mathbf{v}_0(u_0)\|_{B^{1}_{\infty,1}}
&\leq \int^t_0\|\pa_\tau u-\mathbf{v}_0(u_0)\|_{B^{1}_{\infty,1}} \dd\tau
\nonumber\\&\leq \int^t_0\|\partial_x(1-\partial_{xx})^{-1}(u-u_0)\|_{B^{1}_{\infty,1}} \dd\tau+\int^t_0\frac{3}{2}\|u\pa_xu-u_0\pa_xu_0\|_{B^{1}_{\infty,1}} \dd\tau
\nonumber\\&\lesssim \int^t_0\|u(\tau)-u_0\|_{B^{1}_{\infty,1}} \dd\tau+\int^t_0\|u(\tau)-u_0\|_{L^\infty} \|u(\tau)\|_{B^{2}_{\infty,1}} \dd\tau
\nonumber\\&\quad \ + \int^t_0\|u(\tau)-u_0\|_{B^{2}_{\infty,1}}  \|u_0\|_{L^\infty}\dd \tau\nonumber\\
&\lesssim \int^t_0\|u(\tau)-u_0\|_{B^{1}_{\infty,1}} \dd\tau+\|u_0\|_{B^{2}_{\infty,1}}\int^t_0\|u(\tau)-u_0\|_{L^\infty}  \dd\tau
\nonumber\\&\quad \ + \|u_0\|_{L^\infty}\int^t_0\|u(\tau)-u_0\|_{B^{2}_{\infty,1}}  \dd \tau.
\end{align}
Plugging \eqref{et1}--\eqref{et3} into \eqref{et4} yields the desired result \eqref{et0}. Thus, we complete the proof of Proposition \ref{pro1}.

Now we prove the non-uniform continuous dependence.

Set $u^n_0=f_n+g_n$ and compare the solution $\mathcal{S}_{t}(u^n_0)$ with $\mathcal{S}_{t}(f_n)$. Obviously,
\bbal
\|u^n_0-f_n\|_{B^{1}_{\infty,1}}=\|g_n\|_{B^{1}_{\infty,1}}\leq C2^{-n},
\end{align*}
which means that
\bbal
\lim_{n\to\infty}\|u^n_0-f_n\|_{B^{1}_{\infty,1}}=0.
\end{align*}
From Lemma \ref{lem4.1}, one has
\bbal
&\|u^n_0,f_n\|_{B^{\sigma}_{\infty,1}}\leq C2^{(\sigma-1)n}\quad \text{for}\quad \sigma\geq1,\\
&\|u^n_0,f_n\|_{L^\infty}\leq C2^{-n},
\end{align*}
which implies
\bbal
\mathcal{E}(u^n_0)+\mathcal{E}(f_n)\leq C.
\end{align*}
Using the facts
\bbal
&\big\|u^n_{0}\pa_xg_n\big\|_{B^{1}_{\infty,1}}\leq C\big\|u^n_0\big\|_{B^{1}_{\infty,1}}\big\|g_n\big\|_{B^{2}_{\infty,1}}\leq C2^{-n},\\
&\big\|\partial_x(1-\partial_{xx})^{-1}(u_0^n-f_n)\big\|_{B^{1}_{\infty,1}}\leq C\big\|g_n\big\|_{B^{1}_{\infty,1}}\leq C2^{-n}.
\end{align*}
we deduce that
\bal\label{yyh}
\big\|\mathcal{S}_{t}(u^n_0)-\mathcal{S}_{t}(f_n)\big\|_{B^{1}_{\infty,1}}\geq&~ t\big\|g_n\pa_xf_n\big\|_{B^{1}_{\infty,1}}-t\big\|u^n_{0}\pa_xg_n,\;\partial_x(1-\partial_{xx})^{-1}(u_0^n-f_n)\big\|_{B^{1}_{\infty,1}}-Ct^{2}-C2^{-n}\nonumber\\
\geq&~ t\big\|g_n\pa_xf_n\big\|_{B^{1}_{\infty,1}}-Ct2^{-n}-Ct^{2}-C2^{-n},
\end{align}
Notice that \eqref{yz3}
\begin{eqnarray*}
      \liminf_{n\rightarrow \infty} \big\|g_n\pa_xf_n\big\|_{B^{1}_{\infty,1}}\gtrsim M_1,
        \end{eqnarray*}
then we deduce from \eqref{yyh} that
\bbal
\liminf_{n\rightarrow \infty}\big\|\mathcal{S}_t(f_n+g_n)-\mathcal{S}_t(f_n)\big\|_{B^{1}_{\infty,1}}\gtrsim t\quad\text{for} \ t \ \text{small enough}.
\end{align*}
This completes the proof of Theorem \ref{the1.2}.

\section*{Data Availability} Data sharing is not applicable to this article as no new data were created or analyzed in this study.

\section*{Conflict of interest}
The authors declare that they have no conflict of interest.

\section*{Acknowledgments}
 Y. Xiao is supported by the National Natural Science Foundation of China under Grant 11901167.

\end{document}